\documentclass{amsart}
\newtheorem{theorem}{Theorem}[section]
\newtheorem{lemma}[theorem]{Lemma}

\newtheorem{proposition}{Proposition}
\theoremstyle{definition}
\newtheorem{definition}[theorem]{Definition}

\theoremstyle{remark}

\numberwithin{equation}{section}

\begin{document}

\title{Kapranov rank vs. tropical rank}

\author{K. H. Kim}
\address{Department of Mathematics, Alabama State University, Montgomery, Alabama 36101-0271, and
Fellow, Korean Academy of Science and Technology}
\email{khkim@alasu.edu}

\author{F. W. Roush}
\address{Department of Mathematics, Alabama State University,
Montgomery, Alabama 36101-0271}
\email{froush@alasu.edu}
\thanks{}


\date{}

\dedicatory{}


\begin{abstract}
 We show that determining Kapranov rank of tropical matrices is 
not only NP-hard over any infinite field but if solving Diophantine equations over the 
rational numbers is undecidable, then determining Kapranov rank over the 
rational numbers is also undecidable.  We prove that Kapranov rank 
of tropical matrices is not bounded in terms of tropical rank, 
answering a question of Develin, Santos, and Sturmfels \cite{[Stur1]}.
\end{abstract}

\maketitle

\section{Introduction}

The tropical semiring $(R \cup \{ \infty \} ,\min (a,b), a+b )$ 
has many applications but is of special interest in relation 
to algebraic geometry in terms of valuations.  Roughly speaking 
a tropical matrix can arise as the matrix of valuations of a 
matrix of polynomials or power series over a field.  In the paper \cite{[Stur1]} 
three natural notions of rank are introduced.

\begin{definition}
 An $n\times n$ tropical matrix $M$ is nonsingular 
if and only ifthe minimum value over permutations $\pi $ of
$$\sum_{i=1}^n m_{i\pi (i)}$$
is realized uniquely.
\end{definition}

If a matrix of power series is singular, then the corresponding 
tropical matrix must be singular in this sense.

\begin{definition}
The tropical rank of a tropical matrix $M$ is the largest
$k$ such that $M$ has a nonsingular $k\times k$ submatrix.
\end{definition}

The relationship of the tropical semiring to algebraic geometry 
starts with a field of power series $K$, where the power series are 
allowed to have general real exponents in $t$ and the coefficients lie in a base field $F$.
The degree map gives a valuation of this ring.  Let $\mathcal{I}$ be an 
ideal in $K[x_1,\dots ,x_d]$, corresponding to an algebraic variety.  
This gives rise to a variety $V(\mathcal{I})$ in $K^d$ consisting of 
(nonzero) $d$-tuples which when substituted for the $x_i$ make 
all elements of the ideal $\mathcal{I}$ equal to $0$.  That gives rise to a 
tropical variety $T(\mathcal{I})$ in Euclidean $d$-space by taking its 
image under the degree mapping.

For purposes of defining rank the following definition is satisfactory,
however \cite{[SS2]} gives a more general concept using Grassmannians
and hyperplane intersections.

\begin{definition}
A tropical linear space is a tropical variety 
$T(\mathcal{I})$ where $\mathcal{I}$ is generated by 
linear forms in the $x_i$ .  Its dimension 
is $d$ minus the number of minimal generators of $\mathcal{I}$.
\end{definition}

\begin{definition}
The Kapranov rank of a tropical matrix $M$ is the 
least dimension of a tropical linear space containing the rows of $M$.
\end{definition}

The Kapranov rank may differ from the tropical rank, but 
the tropical rank is a lower bound for it.  The Kapranov 
rank is closer to the concept one would like to use 
in algebraic geometry but does not have such a simple tropical 
nature.  Moreover it can vary according to the field.

\begin{definition} 
The Barvinok rank of an $n \times n$ 
tropical matrix $M$ is the least $k$ such that $M$ 
is the product of an $n\times k$ tropical matrix and a $k\times n$ tropical matrix.
\end{definition}

The Barvinok rank is an upper bound for the Kapranov rank.  It is the 
most natural concept of rank for Boolean matrices, where it 
also has names such as Boolean rank and Schein rank \cite{[Kim]}.  
The $2$-element Boolean algebra is the subsemiring of 
the tropical semiring consisting of zero and infinity.  However 
if we look at any $2$ elements of the tropical semiring, the 
additive semigroup is isomorphic to the $2$-element 
Boolean algebra, and some properties of matrices 
whose entries lie in this subset will be like those of Boolean matrices.
For basic properties of Boolean matrices, see \cite{[Kim]}.

Determining Kapranov rank over a given field involves 
dealing with nonlinear systems of polynomial equations 
over that field.  Here we will largely deal with the 
case of Kapranov rank of $(0,1)$-matrices.  Develin, Santos, 
and Sturmfels proved the following.

\begin{theorem}
 \cite{[Stur1]}.  The tropical rank of $M$ is at most $r$ if and 
only if $M$ lies in the variety $T(\mathcal{J_r})$ where 
$\mathcal{J_r}$ is the algebraic variety of matrices 
of the same size as $M$ over $K$ of rank $r$, and $T$ denotes tropicalization.
\end{theorem}

That means that there must be some matrix of rank $r$ over 
the field whose image is $M$.  They deduce the corollary 
that the Kapranov rank is the smallest rank of any lifting over $K$.

\section{Kapranov rank of matrices related to projective planes}

\begin{proposition}
 The Kapranov rank over a infinite field $F$ of 
an $n\times n$  $ (0,1)$-matrix $M$ is the least $k$ 
such that there exist two families of $k$-dimensional vectors $v\langle 
s \rangle , w \langle s\rangle$  over $F$ such that $m_{ij}=0$ 
if and only if the inner product
$v\langle i \rangle   w \langle j\rangle^T=0$ where $T$ denotes transpose.

More generally the Kapranov rank over a field $F$ of an
$n\times n$ nonnegative matrix $M$ is the least
$k$ such that there exist two families of vectors
$v\langle 
s \rangle , w \langle s\rangle$  over the ring $K$ of
power series over $F$ involving arbitrary nonnegative powers
of an indeterminate $t$, such that $m_{ij}$ is precisely 
the order in $t$ of the inner product
$v\langle i \rangle   w \langle j\rangle^T=0$ 
and here $T$ denotes transpose.
\end{proposition}

\begin{proof}
All entries of $M$ lie in the ring $K_+$ 
of power series with no negative degrees.  This is a 
valuation ring in a valuated field and is a principal 
ideal domain. Over this ring, modules of $K$-dimension 
$k$ have also $K_+$ dimension $k$, by the theory of 
torsion-free finitely generated modules over a principal 
ideal domain.  Therefore a matrix representing $M$ will 
factor as a product of $n\times k$ and $k\times n$ matrices, $L$ and $R$ 
over $K_+$, one of which gives the module basis for the row space of $M$ and the 
other gives the linear combinations of that basis which
products all row vectors of $M$.  Then the leading terms of $L$ and $R$ will 
give the required vectors if $M$ does have Kapranov rank at most $k$.
In the second case, the converse is immediate.

Conversely (in the former case) 
suppose that the family of vectors exists.  If we add 
generic degree $1$ terms to $L$ and $R$ then the product will 
have degree $1$ wherever the degree is not zero, and we will 
have the required pattern. 
\end{proof}

This result extends to the case of any two-valued tropical 
matrix.  It would be of some interest to know whether if $M$ 
is nonnegative integer-valued, then one can always work with 
power series having nonnegative integer exponents.  

Instead of working with two families of vectors,
it is equivalent to work with a collection of vectors 
and hyperplanes through the origin, with the condition that 
a vector lies in the hyperplane, if we choose the hyperplanes 
orthogonal to the respective vectors in one of the two sets.
That represents some concept of incidence which is natural
for projective planes.

The following result is also essentially in \cite{[Stur1]}.

\begin{proposition}
 Suppose we have a geometry consisting of a finite 
collection of points and lines and their incidence (membership) 
matrix, and suppose there is at most one line through any two 
points and any two lines intersect in at most one point.
  The tropical rank of the incidence matrix of the 
geometry, taken as a tropical $(0,1)$-matrix, is at most 3. This
remains true if we replace the $1$ entries by arbitrary 
positive real numbers.
\end{proposition}

\begin{proof}
Consider the nature of $4\times 4$ submatrices.  By the Hall-Koenig theorem \cite{[H]}
on systems of distinct representatives 
there is at least one system of distinct representatives 
(matching) for the corresponding binary relation, which consists entirely
of zero unless there are configurations of incidence which 
include $1\times 4, 2\times 3, 3\times 2, 4\times 1$ permuted rectangles
where all points are on all lines.  Of these only $1\times 4$ and $4\times 1$
will occur, and we defer consideration of this case.

 So assume the minimum permuted diagonal sum is zero.  
From \cite{[Kim]}
on matrices of permanent 1 (the complementary result
for Boolean matrices), it follows that if any nonnegative matrix
has a unique permuted diagonal sum of zero, then one can permute the rows and columns
so that the zeroes lie entirely on or below the main diagonal, and the main
diagonal is zero. (This is proved by putting the permuted diagonal on the main diagonal
and noting that the graph of the remaining zeroes must be acyclic, hence
orderable). But such a
pattern means that 2 columns have 2 positive entries in common, two points lie on 
two different lines. 

Now consider the case when the top row consists entirely of positive entries, but
no column consists entirely of positive entries.
By the condition of the theorem, all the other rows have at most one positive
entry.  We claim that the minimum permuted diagonal sum equals the minimum
entry in the top row, and that it is not unique.  Since any sum involves a
term from the top row, this is a lower bound on the minimum.  Now delete its
row and column, we have a $3\times 3$ matrix with at most one positive entry 
per row, and no column consisting entirely of positive entries. There cannot be
any blocking matrix, so this has minimum permuted diagonal of zero.  Moreover
it cannot be permuted so that the zeroes are all on or below the diagonal, so
it has at least two permuted diagonals of zero, and is singular.

Finally consider the case where the top row and the first column are entirely
positive.  Then there can be no other positive entries in the matrix.  
The minimum permuted diagonal sum is either equal
to the $(1,1)$ entry or the minimum sum of another entry in the first row and
another entry in the first column.  In either case, by deleting the rows and columns
of those entries, one finds there will be at least two ways to complete the
permuted diagonal with zeroes.  
\end{proof}

For modules over the ring of power series, it will not only be true 
that a submodule of a finitely generated free module will be free,
but that if we have any spanning set for the submodule, some subset
of that spanning set will be a basis (spanning and independent over
the quotient field). One way to see this is to look at the first coordinate,
take a spanning set element of minimum order there, use it as a first
basis element, and use it to clear out that coordinate entirely; the
resulting module is the kernel of a nontrivial mapping the projection
to that coordinate, and has lower rank over the quotient field, and
one can repeat the process.  This is useful in analyzing certain
configurations.

\begin{lemma}
 Let $V$ be the algebraic set (subset of an algebraic
variety) of all inner products from 2 sequences
of $n$ vectors in $c$-space.  Its tropical dimension does not exceed its ordinary
dimension which is at most a constant times $n$.
\end{lemma}

\begin{proof}
This follows from the Bieri-Groves theorem \cite{[BG]}, alternatively from
the more precise results about tropical Grassmannians in \cite{[SS]}. Note that this
set lies in the algebraic variety of $n\times n$ matrices of rank at most $c$; 
their row and column spaces
represent $c$ planes in $n$-space; the entire set is an image of a product
of two sets of dimension $cn$. 
\end{proof}

\begin{theorem}
Under the above hypothesis, 
there are tropical matrices of tropical rank 3 and arbitrarily high
Kapranov rank.  These are obtained by taking the incidence 
matrices of finite projective planes and inserting general positive real numbers for
the 1 entries (incidences of points on lines). 
\end{theorem}

\begin{proof}
By Proposition 2 these matrices have tropical rank 3. The proof that
they have arbitrarily large Kapranov rank is essentially done by counting dimensions.
Such a matrix $M$ of this can have arbitrarily large size $(m^2+m+1)\times (m^2+m+1)$ 
and $(m^2+m+1)(m+1)$ arbitrary positive entries, since each line has $m+1$ points on it.
Here $m$ can be any order of a finite field, a power of a prime.  If the Kapranov
rank is at most $c$ then by Proposition 1 there are two families $P,L$
called points and lines, each consisting of $n=m^2+m+1$ non-zero vectors
in $c$ dimensional space over a power series ring $R$ over the complex numbers, in a
variable $t$, such that the orders in $t$ of the inner products 
of the points with the lines are 
given by matrix $M$.  All inner products in this set have ordinary dimension
at most $2nc$ by choice of each of $c$ coordinates for
each of $n$ points and lines, but the tropical dimension 
is equal to the number of incidences
in the projective plane, asymptotically $n^{3/2}$.  Therefore Lemma 3
gives a contradiction. 
\end{proof}

\section{Diophantine equations}

Hypothesis. Diophantine equations over the rational numbers are not 
decidable by a general algorithm.

This problem is unknown and very deep; Davis, Putnam, and 
Robinsion, and Matijasevitch proved (Hilbert's Tenth Problem) 
that Diophantine equations over the integers are algorithmically 
undecidable.  To relate this to Kapranov rank, we note by the 
above that representing any system of incidence of lines and 
points in a projective plane is a problem of 
deciding if Kapranov rank is at most $3$.  
Then we make use of Hall's method given in \cite{[Hall]} for 
coordinatizing projective planes and defining algebraic 
operations solely in terms of systems of intersections of points and lines.

The method of Hall to coordinatize a projective plane starts by 
choosing $4$ points $X,Y,O,I$ no three collinear, where $O$ is 
assigned coordinates $(0,0$), $I$ is given coordinates $(1,1)$, $ OX$ 
and $OY$ are the $x$ and $y$-axes, and $XY$ is understood as the line at infinity.
The points on $OI$ other than its intersection denoted (1) with 
$XY$ are given coordinates $(b,b)$.  All points $P$ of the plane 
not on $XY$ are given their $x$ and $y$ coordinates by intersecting  
$XP$ (a horizontal line) and $YP$ (a vertical line) with $OI$.  
Addition is defined by setting $y=x+b$ if and only if $(x,y)$ 
is a finite point on the line joining $(1)$ with $(0,b)$.
Multiplication is defined by setting $y=xm$ if and only if 
$(x,y)$ is a finite point on the line joining $O$ and the infinite 
point $(m)$ which is the intersection of the line joining $(0,0)$ 
and $(1,m)$, with $XY$.  These operations agree with standard 
operations of arithmetic for the standard projective plane over 
any field.  We can take a projective 
transformation so that any given finite set of points and lines 
are finite in the since of being and not lying on $XY$.

For further information on Diophantine undecidability over the rational 
numbers, see \cite{[KR1]}.

\begin{theorem}
Determining whether a two-valued tropical matrix has Kapranov rank 3
over an infinite field has precisely the computational complexity of solving a
system of Diophantine equations over that field.
Under the above hypothesis, computing Kapranov rank 
over the rational numbers is algorithmically undecidable.  Computing 
Kapranov rank of an $n\times n$ $(0,1)$-matrix 
over any infinite field is NP-hard; over an algebraically closed field
it is PSPACE-easy.
\end{theorem}

\begin{proof}
For a $(0,1)$-matrix, having Kapranov rank 3 over the rational numbers is equivalent
by Proposition 1 to finding two systems of 3 dimensional vectors over the field
in question whose
zero and nonzero inner products 
represent it, and as remarked above, to finding a system of nonzero
vectors and planes through the origin in 3 space whose incidence relations
represent it.  For a finite system, we can move all the vectors and planes
away from infinity, and represent them in terms of a system of points
and lines in the plane (or projective plane) over that field.

We use Hall's construction to represent a general Diophantine 
equation with positive integer coefficients on each side as a 
question of existence of a system of points and lines with 
suitable intersections and nonintersections.  We choose 4 points 
as above to define the system, with the non-collinearity assumption.  
Then we choose a number of general points sufficient to define the 
number of variables required.  Then we take the lines and intersections 
required to define each side of the equations as a point on $OI$.
Existence of this configuration is then equivalent to solving the 
Diophantine system.  Some incidences and intersections 
have been undefined, but if we take all possibilities for $0$ and 
$1$ on the incidence matrix for them, it cannot be possible to 
algorithmically decide all the cases.

The same argument proves NP-hardness in general, given that solving 
Diophantine equations or systems over a
field, of polynomial size with coefficients and exponents having 
polynomially many digits is NP-hard.  This follows by taking 
equations representing a Boolean satisfiability problem, adding 
equations $x^2=x$ to force the solutions to be $0$ or $1$. 

However in the NP case it is not sufficient to allow all possibilities
for the unknown incidences, since there might be exponentially many.
Instead, we can control them by adding new variables and equations, and
then taking linear combinations of the given equations.  Our addition and
multiplication constructions define certain points and lines and incidences of these.
In particular cases there may be additional incidences of a point and a line. But if
so this will be because of the separate constructions involved in one or two
equations, not the entire system if there are many equations.  Transform the equations
as in \cite{[KR]} to be linear-quadratic, and arrange that the nth variable instead of
being $0,1$ has $2^n+1$ added to it so that it is either $2^{n+2}+1$ or $2^{n+2}+2$.

  Assume first the field has transcendentals, indeterminate or generic
elements. Choose all equations consist of one side set equal 
to zero. Add 2 new variables, and
2 new equations $E_1,E_2$ in them, which in effect define the new variables as generic
and independent quadratic combinations of the other variables.  
Replace the previous equations by their sums 
with generic linear combinations of the new equations $E_i$; the
combined equations have the same quadratic form as the $E_i$.  
Do the addition and multiplication constructions so as to 
compute in some fixed order all the new
coefficients in each equation and their products with products of at most 2 previous
variables and then add them successively, and finally add in the linear
terms in the 2 new variables.  The coordinates of points and coefficients of line
equations will represent partial steps in the computation, or separate variables
or certain fixed constants or transcendentals for the coefficients.  It will not
matter for our proof exactly what the incidences in each equation or 
pair of equations are, just that they can be computed
in polynomial time and that they will not depend on which 
variables have one of the two values versus which
variables have the other.  

Note that we do not have to construct the 
numerical values of the variables--their existence alone is what is needed to see whether
this configuration is possible and to see whether the Diophantine equations are solvable.
We also do not have to construct the transcendentals, their existence is enough, and that
the linear combinations which we take are nonzero so that we can recover solutions of the
original equations from solutions of the new equations.  However we must construct the 
coefficients as polynomial combinations of the transcendentals
in them, and all integers in them starting with 1.  In the next two paragraphs we argue
that for 2 different equations, or at two different terms for same equation,
the partial steps in the computation will involve different transcendentals, 
and we can have no new incidence 
among them; for different partial computations
of the same term in the same equation, 
there are only the natural incidences, not ones
depending on values of the variables; the incidences involving variables 
will be only that the same variable has the same
value in two different equations.

The addition construction can be arranged to 
construct in turn, using some parallelism, from $(a,a),(b,b)$,
 $x=a,(a,0),x=y+a,y=b,(a+b,b),x=a+b,
(a+b,a+b)$. Thus it adds the following varying 
finite points in order to add
$a,b$: $(a,a),(b,b),(a,0),(a+b,b)$, and the following finite lines: $x=a,x=y+a,y=b,x=a+b$.
The multiplication construction can be arranged to  
construct in turn $x=1,y=a,(1,a),y=ax,x=b,(b,ab),y=ab,(ab,ab)$.  Thus it adds
the following finite points in order
to multiply $a,b$:$(a,a),(b,b),(0,0),(1,1),(1,0),(1,a),(b,ab),
(ab,ab)$ and these lines: $x=1,y=a,y=ax,x=b,y=ab$.   We can choose
in these that $a$ is a partial computation carried from before, and $b$ is a coefficient
or variable, or the integers -1,1 or 2. We do computations of terms
in order, with the new variables in $E_i$ last; when each term is complete, add
it to the previous partial sum.  Consider a term as some 
$(\sum_{r,s} n_{rs}t_rt_s)x_ix_j$ where $t_r,t_s$ are transcendentals
(there may be only one, and $x_i,x_j$ are variables (there may be one or none), and
$n_{rs}$ are integers.
Lowest powers of $2$ are added first in computing binary integers. For subtraction, we multiply
a coefficient, after computing it, by a fixed $-1$ (a fixed diagram verifies its
sum with 1 is 0).  First the transcendentals $t_i$ are multiplied, and among those, first transcendentals from
the $E_i$, then the multiples by them of the powers of 2 in the binary expansions
of $n_{rs}$, then those are added, then once the entire polynomial
coefficient is computed, it is multiplied times $x_i,x_j$, over which we have little
control beyond knowing their two possible values.

The lines above are either vertical or horizontal or at a 45 degree
angle with intercepts $0,1,a,b,a+b,ab$, or have the form $y=ax$.  The points
lie on pairs of these lines and their coordinates are all $0,1,a,b,a+b,ab$.
This means for incidences, in addition to partial computations being 
equal we need only to look 
at ratios $(a+b)/b$ when a sum is being computed, for incidences on some $y=a_1x$ 
and differences $1-a,b-ab$ when a product is being computed, for incidences on
some $x=y+a$; the second coordinate is always
the more advanced in a product and the less advanced in a sum, 
in the partial computation.
Every line or point
involving an $a$ involves a new transcendental or linear combination 
of transcendentals in a new term: these partial computations will
not be the same as those for another equation or another term in the same 
equation, nor can they equal a variable.  
This is also true for
$a-1,ab-b$. As things are arranged, also new transcendentals will
not cancel from $(a+b)/b$ and will not give the value
of any variable; all the more this cannot happen
when we are adding a previous collection of terms.  
Such $a+b$ cannot occur in the stages of multiplying 
in the variables, and cancellation will not occur in adding in a completed term.
In the process of multiplying out a particular term within itself there cannot be
any unpredictable incidences involving the variables, since multiplying them is done
last.  This shows that all incidences are computable in polynomial time from a 
knowledge of the equations, without knowing the values of the variables.

 Solvability of the new system
means precisely realizability in terms of a system of points and lines over the field, hence
Kapranov rank at most 3 of the corresponding $(0,1)$-matrix.
Hence if there were a P algorithm for Kapranov rank 3, 
then there would be a P algorithm for solving a system of Boolean or Diophantine 
equations in terms of 0,1 which is an NP-complete problem.

If we have a non-algebraically closed but infinite field, then a polynomial number of 
trials of substitutions of field elements for transcendentals one by one
 will result in some case where each of the polynomially many intersections
remain distinct.

Canny's result \cite{[Cann]} implies the last statement 
with proposition 1, that is, he proved that solving a system of polynomial
equations over an algebraically closed field is a problem
lying within PSPACE. 
\end{proof}

This extends to the case in which the matrix consists of 
integers of at most polynomial size.  For larger sizes 
we must solve equations in terms of power series which are zero
over large ranges, and the result is less clear.  If it is a
question of extending linearly after a very few stages, then
it can be done.

It has previously been proved that determining Barvinok rank \cite{[Cela]}
and tropical rank \cite{[KR]} are in general NP-hard problems.  
It would also be interesting to know what the general
form of a nonnegative matrix of small tropical rank is, when it is not given directly
by the pattern of zeroes.

It would also be interesting to know whether the Kapranov ranks over the 
complex numbers of a sequence
projective planes over finite fields represented 
as tropical $(0,1)$-matrices tends to infinity.
If so then there can be no purely tropical definition of a lower bound
on Kapranov rank by which Kapranov rank over every field can be bounded, for
the Kapranov ranks over different fields will not be mutually bounded.
However it appears that there will be a tropically computable lower bound 
which can deal with dimensionality phenomenon of Theorem 4 and in that
sense improves on Kapranov rank.

This is related to the topic of realizable matroids.

\end{document}